\documentclass[a4,11pt]{article}

\usepackage{amsfonts}
\usepackage{amsmath}
\usepackage{amssymb}

\usepackage{amsthm}

\theoremstyle{plain}
 \newtheorem{theorem}{Theorem}[section]
 
 \newtheorem*{theorem*}{Theorem}
 \newtheorem{proposition}[theorem]{Proposition}
 \newtheorem{lemma}[theorem]{Lemma}

\theoremstyle{definition}
 
\theoremstyle{remark}
 \newtheorem{remark}[theorem]{Remark}
 
\numberwithin{equation}{section}

\newcommand{\ts}{\theta^*}
\title{On spherical designs obtained from $Q$-polynomial association schemes}

\author{Sho Suda\\
{\small Division of Mathematics, Graduated School of Information Sciences, Tohoku University,} \\
{\small 6-3-09 Aramaki-Aza-Aoba, Aoba-ku, Sendai 980-8579, Japan}\\
{\small suda@ims.is.tohoku.ac.jp}}

\date{\today}

\begin{document}
\maketitle

%%%%%%%%%%%%%%%%%%%%%%%%%%%%%%%%%%%%%%%%%%%%%%%%%%%%%%%%%%%%%%%%%%%%%%%%%%%%%%%%%%%%%%%%%%%%%%%%%%%%%
\begin{abstract}
We characterize that the image of the embedding of the $Q$-polynomial association scheme into eigenspace by primitive idempotent $E_1$ is a spherical $t$-design in terms of the Krein numbers.  
And we show that the strengths of $P$- and $Q$-polynomial schemes as spherical designs are bounded by constant.   
\end{abstract}
\section{Introduction}
 In the study of symmetric association schemes, 
an important technique is embedding of symmetric association schemes into the unit sphere by the primitive idempotent which has no repeated columns.
In 1977, Delsarte-Goethals-Seidel introduced the concept of spherical $t$-designs in unit sphere.
They showed that a spherical $t$-design with degree $s$ which satisfy $t\geq 2s-2$ carries a $Q$-polynomial association scheme.
And in 2006, Bannai-Bannai showed the same conclusion for an antipodal spherical $t$-design with degree $s$ which satisfy $t=2s-3$.
In both cases, the original spherical design is obtained by embedding of carried association scheme into sphere by certain primitive idempotent.

We consider the spherical designs obtained by embedding of $Q$-polynomial association schemes into first eigenspace.
Their spherical designs are always $2$-designs, 
and it is well known that their spherical designs are $3$-designs if and only if $a_1^*=0$.  
Munemasa in \cite{M} showed the criterion that their spherical designs become $t$-design for $t=4,5$ in terms of Krein numbers of $Q$-polynomial association schemes.
In this paper, we show the criterion in the same situation for any $t$ in terms of Krein numbers of $Q$-polynomial association schemes.
Applying this theorem, we show that the strengths of $P$- and $Q$-polynomial association schemes as spherical designs are at most $8$.
This is a dual theorem that the girths of $P$- and $Q$-polynomial schemes are at most $6$ in \cite[Theorem~8.3.6]{BCN} and in \cite[Corollary~30]{L}.
%%%%%%%%%%%%%%%%%%%%%%%%%%%%%%%%%%%%%%%%%%%%%%%%%%%%%%%%%%%%%%%%%%%%%%%%%%%%%%%%%%%%%%%
\section{Preliminaries}
Let $X$ be a finite set and $\mathcal{R}=\{R_0,R_1,\ldots,R_d\}$ be a set of non-empty subsets of $X\times X$.
Let $A_i$ be the adjacency matrix of the graph $(X,R_i)$.
$(X,\mathcal{R})$ is a symmetric association scheme with class $d$ if the following hold:
\begin{enumerate}
\item $A_0 $ is the identity matrix;
\item $\sum_{i=0}^dA_i=J$, where $J$ is the all ones matrix;
\item $A_i^T=A_i$ for $1\leq i\leq d$;
\item $A_iA_j$ is a linear combination of $A_0,A_1,\ldots,A_d$ for $0\leq i,j\leq d$.
\end{enumerate}
The vector space $\mathcal{A}$ spanned by the $A_i$ is an algebra.
$\mathcal{A}$ is called the Bose-Mesner algebra of $(X,\mathcal{R})$.
Since $\mathcal{A}$ is commutative, there exists primitive idempotents $\{E_0,E_1,\ldots,E_d\}$ where $E_0=\frac{1}{|X|}J$.
Since $\mathcal{A}$ is closed under the entry-wise product, we define the Krein parameters $q_{i,j}^k$ as follows: $E_i\circ E_j=\frac{1}{|X|}\sum\nolimits_{k=0}^d q_{i,j}^k E_k$.
$(X,\mathcal{R})$ is $Q$-polynomial (or cometric) with respect to the ordering $\{E_0,E_1,\ldots,E_d\}$ if the following hold:
$q_{i,j}^k=0$ if $i+j>k$ and  $q_{i,j}^k>0$ if $i+j=k$.
If $(X,\mathcal{R})$ is $Q$-polynomial, we define $a_i^*=q_{1,i}^i$, $b_i^*=q_{1,i+1}^i$, and $c_i^*=q_{1,i-1}^i$. 
It is easy to see that $a_i^*+b_i^*+c_i^*=m$ for all $0\leq i\leq d$ and $a_0^*=c_0^*=b_d^*=0$, $c_1^*=1$.

For a positive integer $t$, a finite non-empty set $X$ in the unit sphere $S^{m-1}$ is called 
a spherical $t$-design in $S^{m-1}$ if the following condition is satisfied:
$$\frac{1}{|X|}\sum\limits_{x \in X}f(x)=\frac{1}{|S^{m-1}|}\int\nolimits_{S^{m-1}}f(x)d\sigma(x)$$
for all polynomials $f(x)=f(x_1,\dots,x_m)$ of degree not exceeding $t$.
Here $|S^{m-1}|$ denotes the volume of the sphere $S^{m-1}$.
The following criterion for $t$-designs is well known.
\begin{lemma}\label{criterion design}{\upshape \cite[Corollary 1]{Si}}
Let $X$ be a finite set in $S^{m-1}$. Then the following inequalities hold for all $i\in \mathbb{N}$:
\begin{align}\label{sidel}
\frac{1}{|X|^2}\sum\limits_{x,y\in X}\langle x,y\rangle^i\geq 
\begin{cases}
\frac{(i-1)!!(m-2)!!}{(m+i-2)!!} & \text{if}\ i \text{ is even},\\
0 & \text{if}\ i\text{ is odd}.
\end{cases}
\end{align}
Moreover equalities hold in (\ref{sidel}) for $1\leq i\leq t$ if and only if $X$ is a spherical $t$-design in $S^{m-1}$.  
\end{lemma}

Let $\alpha=(a_i^*)_{0\leq i\leq d}$, $\beta=(b_i^*)_{0\leq i\leq d}$, $\gamma=(c_i^*)_{0\leq i\leq d}$ 
be sequences of non-negative real numbers satisfying 
\begin{align}\label{sequence}
\prod_{i=0}^{d-1} b_i^*c_{i+1}^*\neq0 \text{ and } a_i^*+b_i^*+c_i^*=m  \text{ for all } 0\leq i\leq d \text{ and } a_0^*=c_0^*=b_d^*=0,\ c_1^*=1.
\end{align}   
We define the $(d+1)$-Catalan matrix $F=F^{\alpha,\beta,\gamma}=(f_{n,k})$ of size $d+1$ where $0\leq n,k\leq d$ by the recurrence
\begin{align}
&f_{0,0}=1, \quad f_{0,k}=0,\label{recursion;1}\\
&f_{n,k}=c_k^*f_{n-1,k-1}+a_k^*f_{n-1,k}+b_k^*f_{n-1,k+1}\label{recursion;2}.
\end{align}
Moreover we can define $f_{n,k}$ for $(n,k)\in \{(x,y)\mid d+1\leq n\leq 2d, 0\leq n+k\leq 2d\}$ by the recurrence (\ref{recursion;2}). 
The numbers $B_n=f_{n,0}$ for $0\leq n\leq 2d$ are said to be the Catalan numbers associated with $\alpha$, $\beta$, $\gamma$.

Then $F^{\alpha,\beta,\gamma}$ is a lower triangle matrix and it is easy to see that for $1\leq n\leq d$
\begin{align*}
f_{n,n}=c_1^*\cdots c_n^*, \quad f_{n,n-1}=c_1^*\cdots c_{n-1}^*(a_1^*+\cdots+a_{n-1}^*).
\end{align*}
Therefore we have for $1\leq n\leq d$
\begin{align}\label{eq;1}c_n^*=\frac{f_{n,n}}{f_{n-1,n-1}},\quad a_n^*=\frac{f_{n+1,n}}{f_{n,n}}-\frac{f_{n,n-1}}{f_{n-1,n-1}}.\end{align}

We consider the weighted directed graph $G=(V,E,w)$ associated with $\alpha$, $\beta$, $\gamma$ of non-negative real numbers satisfying (\ref{sequence}) 
where 
\begin{align*}
&V=\{(x,y)\in \mathbb{Z}\times\mathbb{Z}\mid 0\leq y\leq x, x+y\leq2d \},\\ 
&E=\{((x_1,y_1),(x_2,y_2))\in V\times V\mid x_1+1=x_2,|y_1-y_2|\leq1\},\\ 
&w:E\rightarrow \{a_k^*,b_k^*,c_k^*\mid 0\leq k\leq d\}; e\mapsto \begin{cases}
     b_k^* & \text{ if }\ e=((n,k+1),(n+1,k)),\\
     a_k^* & \text{ if }\ e=((n,k),(n+1,k)),\\
     c_k^* & \text{ if }\ e=((n,k-1),(n+1,k)).
\end{cases}
\end{align*}
Let $P$ be a path from $(0,0)$ to $(n,k)\in V$,  
we define the weight of $P$ by 
$$w(P)=\prod_{e\in P}w(e).$$
We define the set of paths from $(0,0)$ to $(n,k)$ by $\mathcal{P}_{n,k}$, the set of the paths from $(0,0)$ to $(n,k)$ via $(\alpha,\beta)$ by $\mathcal{P}_{n,k}^{\alpha,\beta}$, 
and the set of the paths from $(0,0)$ to $(n,k)$ via both $(\alpha,\beta)$ and $(\alpha+1,\beta)$ by $\tilde{\mathcal{P}}_{n,k}^{\alpha,\beta}$.
Then by (\ref{recursion;1}) and (\ref{recursion;2}) we obtain
\begin{align}\label{wp}
f_{n,k}=\sum_{P\in \mathcal{P}_{n,k}}w(P)
\end{align}
for $(n,k)\in V$.
\begin{proposition}\label{p;3}
Let $F=F^{\alpha,\beta,\gamma}$ be a Catalan matrix.
Assume $1\leq t\leq 2d$.
One of the following uniquely determines the others;
\begin{enumerate}
\item the Catalan subsequence $(B_i)_{1\leq i\leq t}$,
\item the submatrix $(f_{n,k})$ where  $(n,k)\in \{(x,y)\in \mathbb{Z}\times\mathbb{Z}\mid 0\leq y\leq x, x+y\leq t \}$,
\item the subsequences $(a_i^*)$ for $0\leq i\leq \lfloor (t-1)/2\rfloor$ and $(c_j^*)$ for $1\leq j\leq \lceil (t-1)/2\rceil$.
\end{enumerate}
\end{proposition}
\begin{proof}
(3) uniquely determines (2) by (\ref{wp}), and (2) uniquely determines (1) by setting $k=0$.
We prove that (1) uniquely determines (3) by induction on $t$.
For $t=1$, there is nothing to prove since $a_0^*=0$.
Let us suppose $t\geq 2$ and that the assertion has been proved for $t-1$.

If $t$ is even, that is $t=2m$, then we have
\begin{align}\label{eq;4}
B_{2m}&=\sum\limits_{P\in \mathcal{P}_{2m,0}}w(P) \displaybreak[0] \notag\\
&=\sum\limits_{k=0}^{m-1}\sum\limits_{P\in \tilde{\mathcal{P}}_{2m,0}^{2m-1-k,k}}w(P)+\sum\limits_{P\in \mathcal{P}_{2m,0}^{m,m}}w(P) \displaybreak[0] \notag\\
&=\sum_{k=0}^{m-1}f_{2m-1-k,k}a_k^*b_{k-1}^*\cdots b_1^*+f_{m-1,m-1}c_m^*b_{m-1}^*\cdots b_1^*.
\end{align} 
By the induction hypothesis, weights $a_1^*,\ldots,a_{m-1}^*$ and $c_1^*,\ldots,c_{m-1}^*$ are uniquely determined by $(B_i)_{1\leq i\leq 2m-1}$.
Since $b_i^*=m-a_i^*-c_i^*$ for any $i$, weights $b_1^*,\ldots,b_{m-1}^*$ are also uniquely determined by $(B_i)_{1\leq i\leq 2m-1}$.
Since (3) uniquely determines (2),  $f_{2m-1,0},f_{2m-2,1},\ldots,f_{m,m-1},f_{m-1.m-1}$ are also uniquely determined by $(B_i)_{1\leq i\leq 2m-1}$.
Therefore it follows from (\ref{eq;4}) that $c_m^*$ is uniquely determined by $(B_i)_{1\leq i\leq 2m}$. 

If $t$ is odd, that is $t=2m+1$, then we have
\begin{align}
B_{2m+1}&=\sum\limits_{P\in \mathcal{P}_{2m+1,0}}w(P)\notag \displaybreak[0]\\ 
&=\sum\limits_{k=0}^{m-1}\sum\limits_{P\in \tilde{\mathcal{P}}_{2m+1,0}^{2m-k,k}}w(P) \notag \displaybreak[0]\\
&=\sum_{k=0}^mf_{2m-k,k}a_k^*b_{k-1}^*\cdots b_1^*\notag \displaybreak[0]\\ 
&=\sum_{k=0}^{m-1}f_{2m-k,k}a_k^*b_{k-1}^*\cdots b_1^*+f_{m,m}a_m^*b_{m-1}^*\cdots b_1^*.\label{eq;5}
\end{align} 
By the induction hypothesis, weights $a_1^*,\ldots,a_{m-1}^*$ and $c_1^*,\ldots,c_m^*$ are uniquely determined by $(B_i)_{1\leq i\leq 2m}$.
Since $b_i^*=m-a_i^*-c_i^*$ for any $i$, weights $b_1^*,\ldots,b_{m-1}^*$ are also uniquely determined by $(B_i)_{1\leq i\leq 2m}$.
Since (3) uniquely determines (2),  $f_{2m,0},f_{2m-1,1},\ldots,f_{m,m}$ are also uniquely determined by $(B_i)_{1\leq i\leq 2m}$.
Therefore it follows from (\ref{eq;5}) that $a_m^*$ is uniquely determined by $(B_i)_{1\leq i\leq 2m+1}$.
\end{proof}

We define the polynomials $v_i^*(x)$ with degree $i$ as follows;
\begin{align}\label{polynomial}
v_0^*=1,\quad v_1^*(x)=x, \quad v_{k+1}^*(x)=\frac{1}{c_{k+1}^*}(xv_k^*(x)-a_k^*v_k^*(x)-b_{k-1}^*v_{k-1}^*(x)),
\end{align}
for $1\leq k\leq d-1$.
Then we can easily find that 
\begin{align}\label{expansion}
x^n=\sum_{k=0}^nf_{n,k}v_k^*(x)
\end{align}  
for $0\leq n\leq d$.
This means that $f_{n,k}$ appears in the coefficients of $\{x^n\}_{0\leq n\leq d}$ expressed in terms of the polynomials $\{v_k^*(x)\}_{0\leq k\leq d}$.

We define the Gegenbauer polynomials $\{Q_k(x)\}_{k=0}^{2d}$ on $S^{m-1}$ by
\begin{align*}
& Q_0(x)=1,\quad Q_1(x)=mx,\\
& \frac{k+1}{m+2k}Q_{k+1}(x)=xQ_k(x)-\frac{m+k-3}{m+2k-4}Q_{k-1}(x),
\end{align*}
for $1\leq k\leq 2d-1$.
It is easily shown that $t^n=\sum_{k=0}^ng_{n,k}Q_k(x)$ where
$$g_{n,k}=\begin{cases}
     \frac{n!(m-2)!!}{(n-k)!!(m+n+k-2)!!} & \text{ if }\ n\equiv k \pmod{2}, \\
     0 & \text{ if }\ n\not\equiv k \pmod{2}. 
\end{cases}$$
Therefore $t^n=\sum_{k=0}^nh_{n,k}Q_k(\frac{x}{m})$ where
\begin{align*}
h_{n,k}=
\begin{cases}
     \frac{m^nn!(m-2)!!}{(n-k)!!(m+n+k-2)!!} & \text{ if }\ n\equiv k \pmod{2}, \\
     0 & \text{ if }\ n\not\equiv k \pmod{2}. 
\end{cases}
\end{align*}
\begin{proposition}\label{p;4}
Let $F=F^{\alpha,\beta,\gamma}$ be a $(d+1)$-Catalan matrix.
The following are equivalent;
\begin{enumerate} 
\item $f_{n,0}=\begin{cases}
     \frac{m^n(n-1)!!(m-2)!!}{(m+n-2)!!}   & \text{ if }\ n \text{ is even },\\
     0 & \text{ if }\ n \text{ is odd }, 
\end{cases}$
for $1\leq n\leq t$,
\item $f_{n,k}=\begin{cases}
     \frac{m^nn!(m-2)!!}{(n-k)!!(m+n+k-2)!!} & \text{ if }\ n\equiv k \pmod{2}, \\
     0 & \text{ if }\ n\not\equiv k \pmod{2}, 
\end{cases}$
for $(n,k)\in \{(x,y)\in \mathbb{Z}\times\mathbb{Z}\mid 0\leq y\leq x, x+y\leq t \}$,
\item $a_i^*=0$ for $0\leq i\leq \lfloor (t-1)/2\rfloor$ and $c_j^*=\frac{mj}{m+2j-2}$ for $1\leq j\leq \lceil (t-1)/2\rceil$.
\end{enumerate}
\end{proposition}
\begin{proof}
By Propositoin~\ref{p;3}, it suffices to show (2)$\Rightarrow $(1), (3).
(2) implies (1) by setting $k=0$ and (2) implies (3) by (\ref{eq;1}).
\end{proof}
%%%%%%%%%%%%%%%%%%%%%%%%%%%%%%%%%%%%%%%%%%%%%%%%%%%%%%%%%%%%%%%%%%%%%%%%%%%%%%%%%%%%%%%
\section{Spherical designs obtained from $Q$-polynomial association schemes}
Let $(X,\mathcal{R})$ be a symmetric association scheme with class $d$ where $E_1$ has no repeated rows, 
and $\tilde{X}$ the image of the embedding into 
the first eigenspace by $E_1=\frac{1}{|X|}\sum\nolimits_{j=0}^dq_1(j)A_j$.
We set $l_0=0$.
For $i\geq1$, comparing the $(x,y)$-entry with $(x,y)\in R_n$ in 
$$(|X|E_1)^i=\sum_{l_1=0}^d\sum_{l_2=0}^d\cdots \sum_{l_i=0}^dq_{1,l_0}^{l_1}q_{1,l_1}^{l_2}\cdots q_{1,l_{i-1}}^{l_i}|X|E_{l_i}, $$
we have
$${q_1(n)}^i=\sum_{l_1=0}^d\sum_{l_2=0}^d\cdots \sum_{l_i=0}^dq_{1,l_0}^{l_1}q_{1,l_1}^{l_2}\cdots q_{1,l_{i-1}}^{l_i}q_{l_i}(n) .$$
Then 
\begin{align*}
\frac{1}{|X|^2}\sum\limits_{x,y\in \tilde{X}}\langle x,y\rangle^i &=\frac{1}{|X|}\sum_{n=0}^d\left(\frac{q_1(n)}{m}\right)^ik_n \displaybreak[0]\\
&=\frac{1}{|X|m^i}\sum_{l_1=0}^d\sum_{l_2=0}^d\cdots \sum_{l_i=0}^dq_{1,l_0}^{l_1}q_{1,l_1}^{l_2}\cdots q_{1,l_{i-1}}^{l_i}\sum_{n=0}^dq_{l_i}(n)k_n\displaybreak[0]\\
&=\frac{1}{m^i}\sum_{l_1=0}^d\sum_{l_2=0}^d\cdots \sum_{l_{i-1}=0}^dq_{1,l_0}^{l_1}q_{1,l_1}^{l_2}\cdots q_{1,l_{i-1}}^{0}
\end{align*}
Therefore it follows from Lemma~\ref{criterion design} that $\tilde{X}$ is a spherical $t$-design in $S^{m-1}$ if and only if 
\begin{align}\label{moment}
\sum_{l_1=0}^d\sum_{l_2=0}^d\cdots \sum_{l_{i-1}=0}^dq_{1,l_0}^{l_1}q_{1,l_1}^{l_2}\cdots q_{1,l_{i-1}}^{0}=\begin{cases}
\frac{(i-1)!!(m-2)!!}{(m+i-2)!!} & \text{if}\ i \text{ is even},\\
0 & \text{if}\ i\text{ is odd}.
\end{cases}
\end{align}
for any $1\leq i\leq t$.

Moreover assume $(X,\mathcal{R})$ is a $Q$-polynomial with respect to the ordering $\{E_0,E_1,\ldots,E_d\}$. 
We define  $\alpha=(a_i^*)_{0\leq i\leq d}$, $\beta=(b_i^*)_{0\leq i\leq d}$, $\gamma=(c_i^*)_{0\leq i\leq d}$ where $a_i^*,b_i^*,c_i^*$ are Krein numbers of  the $Q$-polynomial scheme $(X,\mathcal{R})$.
Then $\alpha,\beta,\gamma$ satisfy (\ref{sequence}), 
so we have the $(d+1)$-Catalan $F^{\alpha,\beta,\gamma}=(f_{n,k})$ matrix  
and the weighted directed graph $G=(V,E,w)$ associated with $\alpha,\beta,\gamma$.    
The left hand side in (\ref{moment}) is 
\begin{align}\label{weight}
\sum_{(l_1,\ldots, l_{i-1})\in Y_i}q_{1,l_0}^{l_1}q_{1,l_1}^{l_2}\cdots q_{1,l_{i-1}}^{0},
\end{align}
where $Y_i=\{(l_1,\ldots,l_{i-1})\in{\{0,\ldots,d\}}^{i-1}\mid 0\leq l_1,l_{i-1}\leq 1, |l_{k-1}-l_k|\leq 1 \text{ for } 1\leq k\leq i-1\}$.
As the element $(l_1,\ldots,l_{i-1})$ of $Y_i$ corresponds to the path $P=(0,0)\sim(1,l_1)\sim\cdots\sim(i-1,l_{i-1})\sim(i,0)$ from $(0,0)$ to $(i,0)$ in $G$, $q_{1,l_0}^{l_1}q_{1,l_1}^{l_2}\cdots q_{1,l_{i-1}}^{0}$ is equal to $w(P)$.  
Thus (\ref{weight}) is equal to 
$$\sum_{P\in\mathcal{P}_{i,0}}w(P)=f_{i,0}.$$
Therefore $\tilde{X}$ is a spherical $t$-design in $S^{m-1}$ if and only if
\begin{align*}
f_{i,0}=\begin{cases}
     \frac{m^i(i-1)!!(m-2)!!}{(m+i-2)!!}   & \text{ if }\ i \text{ is even },\\
     0 & \text{ if }\ i \text{ is odd }, 
\end{cases}
\end{align*}
for any $1\leq i\leq t$.
By Proposition~\ref{p;4}, we obtain the following Theorem;
\begin{theorem}\label{main} 
Let $(X,\mathcal{R})$ be a $Q$-polynomial association scheme with class $d$, 
and $\tilde{X}$ the image of the embedding into 
the first eigenspace by $E_1=\frac{1}{|X|}\sum\nolimits_{j=0}^d\theta_j^*A_j$.
Then the following are equivalent;
\begin{enumerate}
\item $\tilde{X}$ is a spherical $t$-design in $S^{m-1}$, 
\item $a_i^*=0$ for $0\leq i\leq \lfloor (t-1)/2\rfloor$ and $c_j^*=\frac{mj}{m+2j-2}$ for $0\leq j\leq \lceil (t-1)/2\rceil$.
\end{enumerate}
\end{theorem}
%%%%%%%%%%%%%%%%%%%%%%%%%%%%%%%%%%%%%%%%%%%%%%%%%%%%%%%%%%%%%%%%%%%%%%
\section{$P$- and $Q$-polynomial schemes}
We will show the strengths of the spherical design obtained by embedding $P$- and $Q$-polynomial schemes into the first eigenspace are at most $8$.
For general information about distance-regular graphs, see \cite{BI}{BCN}. 
Let $(X,\mathcal{R})$ be a $Q$-polynomial association scheme, and $\tilde{X}$ the image of the embedding into the first eigenspace by $E_1=\frac{1}{|X|}\sum\nolimits_{j=0}^d\theta_j^*A_j$.
For $z\in X$ and $i\in \{1,\ldots,d\}$, 
$\tilde{X}_i(z)$ will denote the orthogonal projection of $R_i(z)$ to $z^\perp=\{y\in \mathbb{R}^{\theta_0^*-1}\mid \langle y,z\rangle=0$, 
rescaled to lie in $S^{\theta_0^*-2}$ in $z^\perp$.
$\tilde{X}_i(z)$ is called the derived design.
In fact $\tilde{X}$ is a ($t+1-d^*$)-design in $S^{\theta_0^*-2}$ by \cite[,Theorem~8.2]{DGS}, where $d^*=\{j\mid 1\leq j\leq d, \ts_j\neq -\ts_0\}$  

The following lemma is used to prove Lemma~\ref{derived}.
\begin{lemma}\label{as}
Let $(X,\mathcal{R})$ be a 
symmetric association scheme of class $d$.
Let $B_i=(p_{i,j}^k)$ be its $i$-th intersection matrix, and $Q=(q_j(i))$ be the second eigenmatrices of 
$\mathfrak{X}$.
Then 
\[
(Q^tB_i)(h,i)=\frac{k_iq_h(i)^2}{m_h}
\quad(0\leq h,i\leq d).
\]
\end{lemma}
\begin{proof}
See \cite[p.73 (4.2) and Theorem 3.5(i)]{BI}.
\end{proof}
The following Lemma gives a condition of derived designs of embedding of a $Q$-polynomial association scheme into first eigenspace
\begin{lemma}\label{derived}
Let $(X,\mathcal{R})$ be a $Q$-polynomial association scheme, and $\tilde{X}$ the image of the embedding into the first eigenspace by $E_1=\frac{1}{|X|}\sum\nolimits_{j=0}^d\theta_j^*A_j$.
Assume $d\geq5$.
Then the following hold for $i\in\{1,\ldots,d\}$ with $\theta_i^*\neq-\theta_0^*$.
\begin{enumerate}
\item $\tilde{X}_i(z)$ is a $2$-design in $S^{\theta_0^*-2}$ for any $z\in X$ if and only if $a_1^*(\ts_i+1)=0$,
\item Assume $\tilde{X}$ is a $4$-design in $S^{\theta_0^*-1}$. 
Then $\tilde{X}_i(z)$ is a $3$-design in $S^{\theta_0^*-2}$ for any $z\in X$ if and only if $a_2^*\bigl((\ts_0+2){\ts_i}^2+2\ts_0\ts_i-{\ts_0}^2\bigr)=0$,
\item Assume $\tilde{X}$ is a $6$-design in $S^{\theta_0^*-1}$. 
Then $\tilde{X}_i(z)$ is a $4$-design in $S^{\theta_0^*-2}$ for any $z\in X$ if and only if 
$a_3^*\bigl((\ts_0+4)(\ts_0-2){\ts_i}^3-3\ts_0(\ts_0+2){\ts_i}^2+3{\ts_0}^2(\ts_0+2)\ts_i+3{\ts_0}^3\bigr)=0$, 
\item Assume $\tilde{X}$ is a $8$-design in $S^{\theta_0^*-1}$. 
Then $\tilde{X}_i(z)$ is a $5$-design in $S^{\theta_0^*-2}$ for any $z\in X$ if and only if
$a_4^*\bigl((\ts_0+4)(\ts_0+6){\ts_i}^4+4\ts_0(\ts_0+4){\ts_i}^3-6{\ts_0}^2(\ts_0+4){\ts_i}^2-12{\ts_0}^3\ts_i+3{\ts_0}^4 \bigr)=0$. 
\end{enumerate} 
\end{lemma}
\begin{proof}
(1) is already shown in \cite[Lemma~4.2]{Su}.
 
The angle set of $\tilde{X}_i(z)$ consists of 
\[
\frac{\frac{\ts_j}{\ts_0}-\frac{{\ts_i}^2}{{\ts_0}^2}}{1
-(\frac{\ts_i}{\ts_0})^2}
=\frac{\ts_0\ts_j-{\ts_i}^2}{{\ts_0}^2-{\ts_i}^2}\quad
(0\leq j\leq d,\;p_{i,j}^i\neq0),
\]
and $|\tilde{X}_i(z)|=k_i$, where $k_i$ is a valency of $R_i$.
Thus, Lemma~\ref{criterion design} implies that $\tilde{X}_i(z)$ is a $t$-design in $S^{\theta_0^*-2}$ is and only if
\begin{align*}
\frac{1}{k_i}\sum\limits_{j=0}^d\Bigl(\frac{\ts_0\ts_j-{\ts_i}^2}{{\ts_0}^2-{\ts_i}^2}\Bigr)^hp_{i,j}^i= 
\begin{cases}
\frac{(h-1)!!(\theta_0^*-3)!!}{(\theta_0^*+h-3)!!} & \text{if}\ h \text{ is even},\\
0 & \text{if}\ h\text{ is odd}
\end{cases},
\end{align*}
for $h=1,\ldots,t$.
Since for $h\leq t\leq d$
\begin{align*}
\sum\limits_{j=0}^d{\theta_j^*}^hp_{i,j}^i&=\sum\limits_{j=0}^d\sum\limits_{l=0}^hf_{i,l}v_l^*(\theta_j^*)p_{i,j}^i&&\text{(by (\ref{expansion}))}\\
&=\sum\limits_{l=0}^hf_{h,l}\sum\limits_{j=0}^dv_l^*(\theta_j^*)p_{i,j}^i\\
&=\sum\limits_{l=0}^hf_{h,l}(Q^tB_i)(l,i)\\
&=k_i\sum\limits_{l=0}^h\frac{f_{h,l}v_l^*(\theta_i^*)^2}{v_l^*(\theta_0^*)},&&\text{(by Lemma~\ref{as})}
\end{align*}
we have 
\begin{align}\label{eq;8}
\frac{1}{k_i}\sum\limits_{j=0}^d\Bigl(\frac{\ts_0\ts_j-{\ts_i}^2}{{\ts_0}^2-{\ts_i}^2}\Bigr)^hp_{i,j}^i&=\frac{1}{k_i({\ts_0}^2-{\ts_i}^2)^h}\sum\limits_{j=0}^d(\ts_0\ts_j-{\ts_i}^2)^hp_{i,j}^i\displaybreak[0] \notag\\
&=\frac{1}{k_i({\ts_0}^2-{\ts_i}^2)^h}\sum\limits_{n=0}^h\binom{h}{n}(-{\ts_i}^2)^{h-n}{\ts_0}^n\sum\limits_{j=0}^d{\ts_j}^np_{i,j}^i\displaybreak[0] \notag\\
&=\frac{1}{({\ts_0}^2-{\ts_i}^2)^h}\sum\limits_{n=0}^h\binom{h}{n}(-{\ts_i}^2)^{h-n}{\ts_0}^n\sum\limits_{l=0}^n\frac{f_{n,l}v_l^*(\ts_i)^2}{v_l^*(\ts_0)}.
\end{align}
To prove (2), assume $\tilde{X}$ is a $4$-design i.e., $a_1^*=0$, $c_2^*=\frac{2\ts_0}{\ts_0+2}$ hold, and  
hence $b_1^*=\ts_0-1$ holds.
It follows from Proposition~\ref{p;4}, (\ref{polynomial}) and (\ref{eq;8}) that 
\begin{align*}
\frac{1}{k_i}\sum\limits_{j=0}^d\Bigl(\frac{\ts_0\ts_j-{\ts_i}^2}{{\ts_0}^2-{\ts_i}^2}\Bigr)^3p_{i,j}^i=\frac{a_2^*\bigl((\ts_0+2){\ts_i}^2+2\ts_0\ts_i-{\ts_0}^2\bigr)^2}{(\ts_0-1)(\ts_0-\ts_i)(\ts_0+\ts_i)^3({\ts_0}^2-(\ts_0+2)a_2^*)}.
\end{align*}
Therefore $\tilde{X}_i(z)$ is a $3$-design if and only if $a_2^*\bigl((\ts_0+2){\ts_i}^2+2\ts_0\ts_i-{\ts_0}^2\bigr)=0$.

To prove (3), assume $\tilde{X}$ is a $6$-design i.e., $a_1^*=a_2^*=0$,  $c_2^*=\frac{2\ts_0}{\ts_0+2}$ and $c_3^*=\frac{3\ts_0}{\ts_0+4}$ hold, and hence $b_1^*=\ts_0-1$ and $b_2^*=\frac{{\ts_0}^2}{\ts_0+2}$ hold.
Then $\tilde{X}_i(z)$ is a $3$-design.
It follows from Proposition~\ref{p;4}, (\ref{polynomial}) and (\ref{eq;8}) that 
\begin{align*}
\frac{1}{k_i}&\sum\limits_{j=0}^d\Bigl(\frac{\ts_0\ts_j-{\ts_i}^2}{{\ts_0}^2-{\ts_i}^2}\Bigr)^4p_{i,j}^i-\frac{3}{(\ts_i+1)(\ts_i-1)}\\
&=\frac{a_3^*\bigl((\ts_0+4)(\ts_0-2){\ts_i}^3-3\ts_0(\ts_0+2){\ts_i}^2+3{\ts_0}^2(\ts_0+2)\ts_i+3{\ts_0}^3\bigr)^2}{({\ts_0}^2-1)(\ts_0+2)(\ts_0-\ts_i)^2(\ts_0+\ts_i)^4(\ts_0(\ts_0+1)-\ts_0a_3^*)}.
\end{align*}
Therefore $\tilde{X}_i(z)$ is a $4$-design if and only if $a_3^*\bigl((\ts_0+4)(\ts_0-2){\ts_i}^3-3\ts_0(\ts_0+2){\ts_i}^2+3{\ts_0}^2(\ts_0+2)\ts_i+3{\ts_0}^3\bigr)=0$.

To prove (4), assume $\tilde{X}$ is a $8$-design i.e., $a_1^*=a_2^*=a_3^*=0$,  $c_2^*=\frac{2\ts_0}{\ts_0+2}$, $c_3^*=\frac{3\ts_0}{\ts_0+4}$ and $c_4^*=\frac{4\ts_0}{\ts_0+6}$ hold, and hence $b_1^*=\ts_0-1$, $b_2^*=\frac{{\ts_0}^2}{\ts_0+2}$ and $b_3^*=\frac{\ts_0(\ts_0+1)}{\ts_0+4}$ hold.
Then $\tilde{X}_i(z)$ is a $4$-design.
It follows from Proposition~\ref{p;4}, (\ref{polynomial}) and (\ref{eq;8}) that 
\begin{align*}
\frac{1}{k_i}&\sum\limits_{j=0}^d\Bigl(\frac{\ts_0\ts_j-{\ts_i}^2}{{\ts_0}^2-{\ts_i}^2}\Bigr)^5p_{i,j}^i\\
&=\frac{a_4^*\bigl((\ts_0+4)(\ts_0+6){\ts_i}^4+4\ts_0(\ts_0+4){\ts_i}^3-6{\ts_0}^2(\ts_0+4){\ts_i}^2-12{\ts_0}^3\ts_i+3{\ts_0}^4 \bigr)^2}{({\ts_0}^2-1)(\ts_0+4)(\ts_0-\ts_i)^3(\ts_0+\ts_i)^5(\ts_0(\ts_0+2)-(\ts_0+6)a_4^*)}.
\end{align*}
Therefore $\tilde{X}_i(z)$ is a $5$-design if and only if $a_4^*\bigl((\ts_0+4)(\ts_0+6){\ts_i}^4+4\ts_0(\ts_0+4){\ts_i}^3-6{\ts_0}^2(\ts_0+4){\ts_i}^2-12{\ts_0}^3\ts_i+3{\ts_0}^4 \bigr)=0$.
\end{proof}
This Lemma implies that for $d\geq5$ and $s\in\{2,3,4,5\}$, the derived designs obtained from a $Q$-polynomial scheme which is a spherical ($2s-1$)-design are spherical $s$-designs. 
Applying this lemma to $P$- and $Q$-polynomial schemes, we obtain the following theorem.
\begin{theorem}
Let $(X,\mathcal{R})$ be a $P$- and $Q$-polynomial scheme with respect to $E_0,E_1,\ldots,E_d$, and 
$\tilde{X}$ the image of the embedding into the first eigenspace by $E_1=\frac{1}{|X|}\sum\nolimits_{j=0}^d\theta_j^*A_j$.
Let $\tilde{X}$ be a $t$-design in $S^{\theta_0^*-2}$.
Then $t\leq8$ holds.
\end{theorem}
\begin{proof}
Assume $t\geq9$.
Since $\tilde{X}$ is a $d$-distance set and $d$-distance sets in $S^{\ts_0-1}$ are at most $2d$-designs, $d\geq5$ holds.
Fix $z\in X$.
Since $\{j\mid R_j\cap(R_1(z)\times R_1(z))\neq \emptyset\} \subset \{0,1,2\}$, $\tilde{X}_1(z)$ is at most $2$-distance set in $S^{\theta_0^*-2}$. 
Hence $\tilde{X}_1(z)$ is at most a $4$-design in $S^{\theta_0^*-2}$.
On the other hand, Lemma~\ref{derived} implies that the derived designs $\tilde{X}_i(z)$ are $5$-designs in $S^{\theta_0^*-2}$, it contradicts. 
\end{proof}
\begin{remark}
(1) Munemasa \cite{M} gives an infinite series that the embedding $P$- and $Q$-polynomial schemes are spherical $5$-designs.   

(2) If a $P$- and $Q$-polynomial scheme $\Gamma$ is a spherical $7$-design, 
then for any $z\in \Gamma$, the local graph $\Gamma_1(z)$ is a tight $4$-design i.e., an extremal Smith graph.
If one can show that there exist no $P$-polynomial schemes whose local graph is an extremal Smith graph, 
the bound is improved i.e., $t\leq6$. 
\end{remark}
%%%%%%%%%%%%%%%%%%%%%%%%%%%%%%%%%%%%%%%%%%%%%%%%%%%%%%%%%%%%%%%%%%%%%%%%%%%%%%%%%%%%
\section*{Acknowledgements}
The author would like to thank Professor Akihiro Munemasa 
for helpful discussions.
This work is supported by Grant-in-Aid for JSPS Fellows.
%%%%%%%%%%%%%%%%%%%%%%%%%%%%%%%%%%%%%%%%%%%%%%%%%%%%%%%%%%%%%%%%%%%%%%%%%%%%%%%%%%%%%%

\end{document}